\documentclass{amsart}
\usepackage{graphics}
\usepackage{amsfonts,amsmath}

\begin{document}

%
%
%
 \newtheorem{thm}{Theorem}[section]
 \newtheorem{cor}[thm]{Corollary}
 \newtheorem{lem}[thm]{Lemma}
 \newtheorem{prop}[thm]{Proposition}
 \newtheorem{defn}[thm]{Definition}
 \newtheorem{assumption}[thm]{Assumption}
 \newtheorem{rem}[thm]{Remark}
 \newtheorem{ex}{Example}
\numberwithin{equation}{section}
 \def\e{{\rm e}}
\def\F{\mathcal{F}}
\def\R{\mathbb{R}}
\def\T{\mathbf{T}}
\def\N{\mathbb{N}}
\def\K{\mathbf{K}}
\def\Q{\mathbf{Q}}
\def\M{\mathbf{M}}
\def\O{\mathbf{O}}
\def\C{\mathbb{C}}
\def\P{\mathbf{P}}
\def\Z{\mathbb{Z}}
\def\H{\mathcal{H}}
\def\A{\mathbf{A}}
\def\V{\mathbf{V}}
\def\AA{\overline{\mathbf{A}}}
\def\B{\mathbf{B}}
\def\c{\mathbf{C}}
\def\L{\mathbf{L}}
\def\bS{\mathbf{S}}
\def\H{\mathcal{H}}
\def\I{\mathbf{I}}
\def\Y{\mathbf{Y}}
\def\X{\mathbf{X}}
\def\f{\mathbf{f}}
\def\z{\mathbf{z}}
\def\d{\hat{d}}
\def\bx{\mathbf{x}}
\def\y{\mathbf{y}}
\def\w{\mathbf{w}}
\def\b{\mathcal{B}}
\def\c{\mathbf{c}}
\def\a{\mathbf{a}}
\def\u{\mathbf{u}}
\def\s{\mathcal{S}}
\def\cc{\mathcal{C}}
\def\co{{\rm co}\,}
\def\vol{{\rm vol}\,}
\def\om{\mathbf{\Omega}}

\title[parametric polynomial optimization]{A "Joint+marginal" approach to parametric polynomial optimization}

\author{Jean B. Lasserre}
\address{LAAS-CNRS and Institute of Mathematics\\
University of Toulouse\\
LAAS, 7 avenue du Colonel Roche\\
31077 Toulouse C\'edex 4,France}
\email{lasserre@laas.fr}
\date{}

\begin{abstract}
Given a compact parameter set $\Y\subset\R^p$, 
we consider polynomial optimization problems $(\P_\y$) on $\R^n$ whose description depends on 
the parameter $\y\in\Y$.
We assume that one can compute all moments  of
some probability measure $\varphi$ on $\Y$,
absolutely continuous with respect to the Lebesgue measure
(e.g. $\Y$ is a box or a simplex and $\varphi$ is uniformly distributed).
We then provide a hierarchy of semidefinite relaxations 
whose associated sequence of optimal solutions converges to the moment vector of a probability measure
that encodes all information about all global optimal solutions $\bx^*(\y)$ of $\P_\y$, as $\y\in\Y$.
In particular, one may approximate as closely as 
desired any polynomial functional of the optimal solutions, like e.g. their $\varphi$-mean. In addition,
using this knowledge on moments, the measurable function $\y\mapsto x^*_k(\y)$ of the $k$-th coordinate of optimal solutions, 
can be estimated, e.g. by maximum entropy methods. Also, for a boolean variable $x_k$,
one may approximate as closely as desired its persistency $\varphi(\{\y:x^*_k(\y)=1\}$, i.e.
the probability that in an optimal solution $\bx^*(\y)$, the coordinate $x^*_k(\y)$
takes the value $1$.
At last but not least, from an optimal solution of the dual semidefinite relaxations,
one provides a sequence of polynomial (resp. piecewise polynomial) lower approximations 
with $L_1(\varphi)$ (resp. almost uniform) convergence to the optimal value function.
\end{abstract}

\keywords{Parametric and polynomial optimization; semidefinite relaxations}

\subjclass{65 D15, 65 K05, 46 N10, 90 C22}

\maketitle
\section{Introduction}~

Roughly speaking, given a set parameters $\Y$ and an optimization problem 
whose description depends on $\y\in\Y$ (call it $\P_\y$),
{\it parametric optimization} is concerned with 
the behavior and properties of the optimal value as well as primal (and possibly dual)
optimal solutions of $\P_\y$, when $\y$ varies in $\Y$. This a quite challenging problem
and in general one may obtain information locally around
some nominal value $\y_0$ of the parameter. There is a vast and rich literature on the topic and 
for a detailed treatment, the interested reader is referred to e.g. Bonnans and Shapiro
\cite{bonnans} and the many references therein.  Sometimes,
in the context of optimization with data uncertainty, some probability distribution $\varphi$ on the parameter set $\Y$ is 
available and in this context one is also interested in e.g. the distribution of the optimal value, optimal solutions, all viewed as random variables. 
In particular, for discrete optimization problems where cost coefficients are random variables with joint distribution $\varphi$,
some bounds on the expected optimal value have been obtained. More recently Natarajan et al. \cite{natarajan} extended the earlier work in
\cite{persis} to even provide a convex optimization problem for computing the so-called {\em persistency} values\footnote{Given 
a $0-1$ optimization problem $\max\{\c'\bx\,:\bx\in\mathcal{X}\cap\{0,1\}^n\}$ and a distribution $\varphi$ on $\c$,
the persistency value of the variable $x_i$ is ${\rm Prob}_\varphi(x_i ^*=1)$ at an optimal solution $\bx^*(\c)=(x^*_i)$.} of (discrete) variables, for a particular 
distribution $\varphi^*$ in a certain set $\Theta$ of distributions. However, this convex formulation requires knowledge of the convex hull of a discrete set and approximations are needed. 
The approach is nicely illustrated on a discrete choice problem and a stochastic knapsack problem.
For more details on persistency in discrete optimization, the interested reader is referred to \cite{natarajan} and the references therein.

In the context of polynomial equations whose coefficients 
are themselves polynomials of some parameter $\y\in\Y$, some specific "parametric"
methods exist. For instance, one may compute symbolically once and for all, what is called a {\it comprehensive} Gr\"obner basis,
i.e., a fixed basis that is a Gr\"obner basis for all $\y\in\Y$; see Weispfenning \cite{weis} and more recently Rostalski \cite{rostalski} for more details.
Then when needed, one may compute the solutions 
for a specific value of the parameter $\y$, e.g.  by the eigenvalue method of M\"oller and Stetter \cite{moller,stetter}. However, one still needs to
apply the latter method for each value of the prameter $\y$. A similar two-step approach is also proposed for homotopy
(instead of Gr\"obner bases) methods in \cite{rostalski}.\\

The purpose of this paper is to show that
in one restricts to the case of {\it polynomial} parametric optimization then all information about
the optimal value and optimal solutions
can be obtained, or at least, approximated as closely as desired.

\subsection*{Contribution}
We here restrict our attention to parametric polynomial optimization, that is, when
$\P_\y$ is described by {\it polynomial} equality and inequality constraints
on both the parameter vector $\y$ and the optimization variables $\bx$. Moreover, the set
$\Y$ is restricted to be a compact basic semi-algebraic set of $\R^p$, and preferably a set 
sufficiently simple so that one may obtain the 
moments of some probability measure on $\Y$, absolutely continuous with respect to the Lebesgue
measure. For instance if
$\Y$ is a simple set (like a simplex, a box) one may choose $\varphi$ to be the probability measure uniformly distributed on $\Y$; typical $\Y$ candidates are polyhedra. Or sometimes, in the context of optimization with data uncertainty, 
$\varphi$ is already specified. We also suppose that 
$\P_\y$ has a unique optimal solution for almost all values of the parameter $\y\in\Y$. 
In this specific context we are going to show
that one may get insightful information on the set of all global optimal 
solutions of $\P_\y$, via what we call a {\it "Joint+marginal"} approach. Our contribution is as follows:\\

(a) Call $J(\y)$ (resp. $\X^*_\y\in\R^n$) the optimal value (resp. the set of optimal solutions) of  $\P_\y$ for the value 
$\y\in\Y$ of the parameter.
We first define an infinite-dimensional optimization problem $\P$ whose optimal value is 
exactly $\rho=
\int_\Y J(\y)d\varphi(\y)$. Any optimal solution of $\P_\y$ is a probability measure $\mu^*$ on $\R^n\times\R^p$ with marginal $\varphi$ on $\R^p$.
It turns out that $\mu^*$ encodes all information on the optimal solutions $\X^*_\y$, $\y\in\Y$.
Whence the name {\it "Joint+marginal"} as $\mu^*$ is a {\it joint} distribution of $\bx$ and $\y$,
and $\varphi$ is the {\it marginal} of $\mu^*$ on $\R^p$.

(b) Next, we provide a hierarchy of semidefinite relaxations of $\P$ with associated sequence of optimal values $(\rho_i)_i$, in the spirit of the hierarchy defined in \cite{lassiopt}.
An optimal solution of the $i$-th semidefinite relaxation is a sequence 
$\z^i=(z^i_{\alpha\beta})$ indexed in the monomial basis $(\bx^\alpha\y^\beta)$ of the subspace 
$\R[\bx,\y]_i$ of polynomials of degree at most $2i$.
If for almost all $\y\in\Y$, $\P_\y$ has a unique {\it global} optimal solution 
$\bx^*(\y)\in\R^n$, then as $i\to\infty$, 
$\z^i$ converges pointwise to the sequence of moments of $\mu^*$ defined in (a).
In particular,  one obtains the distribution of the optimal solution $\bx^*(\y)$, and
therefore, one may approximate as closely as desired any polynomial functional of the solution 
$\bx^*(\y)$, like e.g. the $\varphi$-mean or variance of $\bx^*(\y)$.

In addition, if the optimization variable $x_k$ is boolean then one may approximate as closely as desired
its {\it persistency} $\varphi(\{\y:x^*_k(\y)=1\}$ (i.e., the probability that $x^*_k(\y)=1$ 
in an optimal solution $\bx^*(\y)$), as well as a a necessary and sufficient condition
for this persistency to be $1$.

(c) Finally, let $e(k)\in\N^n$ be the vector $(\delta_{j=k})_j$. Then as $i\to\infty$, and for every $\beta\in\N^p$,
the sequence $(z^i_{e(k)\beta})$ converges to $z^*_{k\beta}:=\int_\Y \y^\beta g_k(\y)d\varphi(\y)$
for the measurable function $\y\mapsto g_k(\y):=x^*_k(\y)$.
In other words, the sequence $(z^*_{k\beta})_{\beta\in\N^p}$  is the moment sequence of the measure
$d\psi(y):=x_k^*(\y)d\varphi(\y)$ on $\Y$. And so, 
the $k$-th coordinate function $\y\mapsto x_k^*(\y)$ of optimal solutions of $\P_\y$, $\y\in\Y$, 
can be estimated, e.g. by maximum entropy methods. Of course, the latter estimation 
is not pointwise but it still provides
useful information on optimal solutions, e.g. the shape of the function $\y\mapsto x^*_k(\y)$,
especially if the function $x^*_k(\cdot)$ is continuous, as illustrated on some simple examples. 
For instance, for parametric polynomial equations,
one may use this estimation of $\bx^*(\y)$ as an initial point for Newton's method 
for any given value of the parameter $\y$. 

Finally, the computational complexity of the above methodology is roughly the same as the 
moment approach described in \cite{lassiopt} for an optimization problem with $n+p$ variables since we consider the joint distribution of the $n$ variables $\bx$ and the $p$ parameters $\y$.
Hence, the approach is particularly interesting when the number of parameters is small, 
say $1$ or $2$. In addition, in the latter case the max-entropy estimation has been shown to 
be very efficient in several examples in the literature; see e.g. \cite{BMP,tagliani1,tagliani2}.
However, in view of the present status of SDP solvers,
if no sparsity or symmetry is taken into account as proposed in e.g. \cite{lassiopt2}, 
the approach is limited to small to medium size polynomial optimization problems.

But this computational price may not seem that high in view of
the ambitious goal of the approach.
After all, keep in mind that by applying the moment approach
to a single $(n+p)$-variables problem, one obtains information 
on global optimal solutions of an $n$-variables problem that depends on $p$ parameters, that
is, one approximates $n$ {\it functions} of $p$ variables!

\section{A related linear program}

Let $\R[\bx,\y]$ denote the ring of polynomials in the variables $\bx=(x_1,\ldots,x_n)$, and
the variables $\y=(y_1,\ldots,y_p)$, whereas $\R[\bx,\y]_k$ denotes its subspace of polynomials of degree at most $k$.  
Let $\Sigma[\bx,\y]\subset\R[\bx,\y]$ denote the subset of polynomials that are sums of squares (in short s.o.s.). For 
a real symmetric matrix $\A$ the notation $\A\succeq0$ stands for $\A$ is positive semidefinite.

\subsection*{The parametric optimization problem}

Let $\Y\subset\R^p$ be a compact set, called the {\it parameter} set, and 
let $f,h_j\,:\,\R^n\times\R^p\to\,\R$, $j=1,\ldots,m$, be continuous. 
For each $\y\in\Y$, fixed,  consider the following optimization problem: 
\begin{equation}
\label{pb1}
J(\y)\,:=\,\inf_\bx\:\{\,f_\y(\bx)\::\: h_{\y j}(\bx)\,\geq\,0,\:j=1,\ldots,m\,\}
\end{equation}
where the functions $f_\y,h_{\y j}\,:\,\R^n\to\,\R$ are defined via:
\[\left.\begin{array}{lcl}
\bx\mapsto f_\y(\bx)&:=&f(\bx,\y)\\
\bx\mapsto h_{\y j}(\bx)&:=&h_j(\bx,\y),\:j=1,\ldots,m
\end{array}\right\}
\quad\forall\,\bx\in\R^n,\,\forall\,\y\in\R^p.\]
Next, let $\K\subset\R^n\times\R^p$ be the set:
\begin{equation}
\label{set-xy}
\K\,:=\,\{\,(\bx,\y)\::\: \y\in\Y\,;\quad h_j(\bx,\y)\,\geq\,0,\quad j=1,\ldots,m\,\},
\end{equation}
and for each $\y\in\Y$, let 
\begin{equation}
\label{set-x}
\K_\y\,:=\,\{\,\bx\in\R^n\::\: h_{\y j}(\bx)\,\geq\,0,\quad j=1,\ldots,m\,\}.
\end{equation}

The interpretation is as follows: $\Y$ is a set of parameters and for each instance
$\y\in\Y$ of the parameter, one wishes to compute an optimal {\it decision} vector $\bx^*(\y)$ that solves problem (\ref{pb1}). 
Let $\varphi$ be a Borel probability measure on $\Y$, with a positive density
with respect to the Lebesgue measure on $\R^p$. For instance
choose for $\varphi$ the probability measure 
\[\varphi(B)\,:=\,\left(\int_\Y d\y\,\right)^{-1}\displaystyle\int_{\Y\cap B}d\y,\qquad\forall B\in\mathcal{B}(\R^p),\]
uniformly distributed on $\Y$.  Sometimes, e.g. in the context of optimization with data uncertainty, $\varphi$ is already specified.

We will use $\varphi$ (or more precisely, its moments) to 
get information on the distribution of
optimal solutions $\bx^*(\y)$ of $\P_\y$, viewed as random vectors.\\

In the rest of the paper we assume that for every $\y\in\Y$, the set $\K_\y$ in (\ref{set-x})  is nonempty.

\subsection{A related infinite-dimensional linear program}

Let $\M(\K)$ be the set of finite Borel measures on $\K$, and consider the following 
infinite-dimensional linear program $\P$:
\begin{equation}
\label{pb2}
\rho\,:=\,\inf_{\mu\in\M(\K)}\:\{\,\int_\K f\,d\mu \::\: \pi\mu\,=\,\varphi\,\}
\end{equation}
where $\pi\mu$ denotes the marginal of $\mu$ on $\R^p$, that is,
$\pi\mu$ is a probability measure on $\R^p$  defined by
\[\pi\mu(B)\,:=\,\mu(\R^n\times B),\qquad\forall \,B\in\mathcal{B}(\R^p).\]
Notice that $\mu(\K)=1$ for any feasible solution $\mu$ of $\P$.
Indeed, as $\varphi$ is a probability measure and $\pi\mu=\varphi$ one has
$1=\varphi(\Y)=\mu(\R^n\times\R^p)=\mu(\K)$.

Recall that for two Borel spaces $X,Y$,
the graph ${\rm Gr}\psi\subset X\times Y$ of a set-valued mapping $\psi: X\to Y$ is the set
\[{\rm Gr}\,\psi\,:=\,\{(\bx,\y)\::\:\bx\in X\,;\: \y\in \psi(\bx)\:\}.\]
If $\psi$ is measurable then any measurable function $h:X\to Y$ with
$h(\bx)\in\psi(\bx)$ for every $\bx\in X$, is called a (measurable) {\it selector}.

\begin{lem}
\label{selector}
Let both $\Y\subset\R^n$ and $\K$ in (\ref{set-xy}) be compact.
Then the set-valued mapping $\y\mapsto \K_\y$ is Borel-measurable. In addition:

{\rm (a)} The mapping $\y\mapsto J(\y)$ is measurable.

{\rm (b)} There exists  a measurable selector $g\,:\Y\to \K_y$ such that
$J(\y)=f(g(\y),\y)$ for every $\y\in\Y$.
\end{lem}
\begin{proof}
As $\K$ and $\Y$ are both compact, the set valued mapping $\y\mapsto\K_\y\subset\R^n$ is compact-valued. Moreover, the graph of $\K_\y$ is by definition the set $\K$, which is a Borel subset of $\R^n\times\R^p$. Hence, by \cite[Proposition D.4]{ohl}, $\K_\y$ is
a measurable function from $\Y$ to the space of nonempty compact subsets of $\R^n$, topologized by the Hausdorff metric. Next, since $\bx\mapsto f_\y(\bx)$ is continuous for every $\y\in\Y$, (a) and (b) follows from e.g. \cite[Proposition D.5]{ohl}.
\end{proof}

\begin{thm}
\label{th1}
Let both $\Y\subset\R^p$ and $\K$ in (\ref{set-xy}) 
be compact and assume that for every $\y\in\Y$,
the set $\K_\y\subset\R^n$ in (\ref{set-x}) is nonempty.
Let $\P$ be the optimization problem (\ref{pb2}) and let
$\X^*_\y:=\{\bx\in\R^n\,:\,f(\bx,\y)=J(\y)\}$, $\y\in\Y$. Then:

{\rm (a)} $\rho\,=\,\displaystyle\int_\Y J(\y)\,d\varphi(\y)$ and $\P$ has an optimal solution.

{\rm (b)} For every optimal solution 
$\mu^*$ of $\P$, and for almost all $\y\in\Y$, there is a probability measure $\psi^*(d\bx\,\vert\,\y)$ on
$\X^*_\y$ such that:
\begin{equation}
\label{th1-1}
\mu^*(C\times B)\,=\,\int_{B\cap\Y}\psi^*(C\cap\X^*_\y\,\vert\,\y)\,d\varphi(\y),\qquad\forall B\in\mathcal{B}(\R^p),\:C\in\mathcal{B}(\R^n).
\end{equation}
\indent
{\rm (c)} Assume that for almost $\y\in\Y$, the set of minimizers of
$\X^*_\y$ is the singleton 
$\{\bx^*(\y)\}$ for some $\bx^*(\y)\in\K_y$. Then there is a measurable mapping $g:\Y\to\K_\y$ such that
\begin{equation}
\label{th1-3}
g(\y)\,=\,\bx^*(\y)\:\mbox{ for every }\:\y\in\Y\,;\quad\rho\,=\,\int_\Y f(g(\y),\y)\,d\varphi(\y),
\end{equation}
and for every $\alpha\in\N^n$, and $\beta\in\N^p$:
\begin{equation}
\label{th1-4}
\int_\K\bx^\alpha\y^\beta\,d\mu^*(\bx,\y)\,=\,\int_\Y\y^\beta \,g(\y)^\alpha\,d\varphi(\y).
\end{equation}
\end{thm}
\begin{proof}
(a) As $\K$ is compact then so is $\K_\y$ for every $\y\in\Y$. Next, as $\K_\y\neq\emptyset$ for every $\y\in\Y$
and $f$ is continuous, the set $\X^*_\y:=\{\bx\in\R^n\,:\,f(\bx,\y)=J(\y)\}$ is nonempty for every $\y\in\Y$.
Let $\mu$ be any feasible  solution of $\P$ and so by definition,
its marginal on $\R^p$ is just $\varphi$. Since
$\X^*_\y\neq\emptyset,\:\forall\y\in\Y$, one has $f_\y(\bx)\geq J(\y)$ for all $\bx\in\K_\y$ and all 
$\y\in\Y$. So,
$f(\bx,\y)\geq J(\y)$ for all $(\bx,\y)\in\K$ and therefore 
\[\int_\K fd\mu\,\geq\,\int_\K J(\y)\,d\mu\,=\,\int_\Y J(\y)\,d\varphi,\]
which proves that $\rho\geq \displaystyle\int_\Y J(\y)\,d\varphi$.

On the other hand, recall that $\K_\y\neq\emptyset,\:\forall\y\in\Y$.
Consider the set-valued mapping $\y\mapsto \X^*_\y\subset\K_\y$. As
$f$ is continuous and $\K$ is compact, then $\X^*_\y$ is compact-valued. In addition,
as $f_\y$ is continuous, by \cite[D6]{ohl} (or \cite{schal}) there exists a measurable selector $g:\Y\to\X^*_y$ (and so $f(g(\y),\y)=J(\y)$).
Therefore, for every $\y\in\Y$, let $\psi^*_\y$  be the Dirac probability measure 
with support on the singleton $g(\y)\in\X^*_\y$, and let $\mu$ be the probability measure
on $\K$ defined by:
\[\mu(C,B)\,:=\,\int_B{\rm 1}_{C}(g(\y))\,\varphi(d\y),\qquad\forall B\in\mathcal{B}(\R^p),\:C\in\mathcal{B}(\R^n).\]
(The measure $\mu$ is well-defined because $g$ is measurable.)
Then $\mu$ is feasible for $\P$ and
\begin{eqnarray*}
\rho\leq\int_\K f\,d\mu&=&\int_\Y \left[\int_{\K_\y}\,f(\bx,\y)\,d\delta_{g(\y)}\,\right]\,d\varphi(\y)\\
&=&\int_\Y f(g(\y),\y)\,d\varphi(\y)=\int_\Y J(\y)\,d\varphi(\y),\end{eqnarray*}
which shows that $\mu$ is an optimal solution of $\P$ and $\rho=\int_\Y J(\y)d\varphi(\y)$.

(b) 
Let $\mu^*$ be an arbitrary optimal solution of $\P$, hence supported on
$\K_\y\times\Y$. Therefore,
as $\K$ is contained in the cartesian product $\R^p\times\R^n$, the probability measure
$\mu^*$ can be disintegrated as
\[\mu^*(C,B)\,:=\,\int_{B\cap\Y}\,\psi^*(C\cap\K_\y\,\vert\,\y)\,d\varphi(\y),\qquad\forall B\in\mathcal{B}(\R^p),\:C\in\mathcal{B}(\R^n),\]
where for all $\y\in\Y$, $\psi^*(\cdot\,\vert\,\y)$ is a probability measure on $\K_\y$.
(The object $\psi^*(\cdot\vert\cdot)$ is called a stochastic kernel; see e.g. 
\cite[p. 88--89]{dynkin} or \cite[D8]{ohl}.) Hence from (a),
\begin{eqnarray*}
\rho=\int_\Y J(\y)\,d\varphi(\y)&=&\int_\K f(\bx,\y)\,d\mu^*(\bx,\y)\\
&=&\int_\Y \left(\int_{\K_\y}f(\bx,\y)\,\psi^*(d\bx\,\vert\,\y)\right)\,d\varphi(y).\end{eqnarray*}
Therefore, using $f(\bx,\y)\geq J(\y)$ on $\K$, 
\[0=\int_\Y\left(\int_{\K_\y} \underbrace{J(\y)- f(\bx,\y)}_{\leq0}\,\psi^*(d\bx\,\vert\,\y)\right)\,d\varphi(y),\]
which implies $\psi^*(\X^*(\y)\,\vert\,\y)=1$ for almost all $\y\in\Y$.

(c) Let $g\,:\Y\to\K_\y$ be the measurable mapping of Lemma \ref{selector}(b).
As $J(\y)=f(g(\y),\y)$ and $(g(\y),\y)\in\K$ then necessarily $g(\y)\in\X^*_\y$ for every $\y\in\Y$.
Next, let $\mu^*$ be an optimal solution of $\P$, and let $\alpha\in\N^n$, $\beta\in\N^p$. Then 
\begin{eqnarray*}
\int_\K\bx^\alpha\y^\beta\,d\mu^*(\bx,\y)&=&
\int_\Y\y^\beta\left(\int_{\X^*_\y}\bx^\alpha\,\psi^*(d\bx\vert\,\y)\right)\,d\varphi(\y)\\
&=&\int_\Y\y^\beta\,g(\y)^\alpha\,d\varphi(\y),
\end{eqnarray*}
the desired result.
\end{proof}

An optimal solution $\mu^*$ of $\P$ encodes {\it all} information on the 
optimal solutions $\bx^*(\y)$ of $\P_\y$.
For instance, let $\B$ be a given Borel set of $\R^n$. Then from Theorem \ref{th1},
\[{\rm Prob}\,(\bx^*(\y)\in\B)\,=\,\mu^*(\B\times\R^p)\,=\,\int_\Y\psi^*(\B\,\vert\,\y)\,d\varphi(\y),\]
with $\psi^*$ as in Theorem \ref{th1}(b).

Consequently, if one knows an optimal solution 
$\mu^*$ of $\P$ then one may evaluate functionals on the solutions of
$\P_\y$, $\y\in\Y$. That is,  assuming that for almost all $\y\in\Y$,
problem $\P_\y$ has a unique optimal solution $\bx^*(\y)$, and given a measurable 
mapping $h\,:\R^n\to\,\R^q$, one may
evaluate the functional
\[\int_\Y h(\bx^*(\y))\,d\varphi(\y).\]
For instance, with $\bx\mapsto h(\bx):=\bx$ one obtains the {\it mean} vector 
${\rm E}_\varphi(\bx^*(\y)):=\int_\Y\bx^*(\y)d\varphi(\y)$ of
optimal solutions $\bx^*(\y)$, $\y\in\Y$.
\begin{cor}
\label{coro1}
Let both $\Y\subset\R^p$ and $\K$ in (\ref{set-xy}) 
be compact. Assume that for every $\y\in\Y$,
the set $\K_\y\subset\R^n$ in (\ref{set-x}) is nonempty, and
for almost all $\y\in\Y$, the set $\X^*_y:=\{\bx \in\K_\y\::\: J(\y)=f(\bx,\y)\}$ is 
the singleton $\{\bx^*(\y)\}$.  Then for every measurable mapping
$h\,:\R^n\to\,\R^q$,
\begin{equation}
\label{coro1-1}
\int_\Y h(\bx^*(\y))\,d\varphi(\y)\,=\,\int_\K h(\bx)\,d\mu^*(\bx,\y).
\end{equation}
where $\mu^*$ is an optimal solution of $\P$. 
\end{cor}
\begin{proof}
By Theorem \ref{th1}(c)
\[\int_\K h(\bx)\,d\mu^*(\bx,\y)\,=\,\int_\Y\left[\int_{\X^*_\y}h(\bx)\psi^*(d\bx\,\vert\,\y)\right]\,d\varphi(\y)
\,=\,\int_\Y h(\bx^*(\y))\,d\varphi(\y).\]
\end{proof}

\subsection{Duality}

Consider the following infinite-dimensional linear program $\P^*$:
\begin{equation}
\label{dual-lp}
\begin{array}{ll}
\rho^*\,:=\,\displaystyle\sup_{p\in\R[\y]}&\displaystyle\int_\Y p\,d\varphi\\
&f(\bx,\y)-p(\y)\,\geq\,0\quad\forall (\bx,\y)\in\K.\end{array}
\end{equation}
Then  $\P^*$ is a {\em dual} of $\P$.

\begin{lem}
\label{lem-dual}
Let both $\Y\subset\R^p$ and $\K$ in (\ref{set-xy}) 
be compact and let $\P$ and $\P^*$ be as in (\ref{pb2}) and (\ref{dual-lp}) respectively.
Then there is no duality gap, i.e., $\rho=\rho^*$.
\end{lem}
\begin{proof}
For a topological space $\mathcal{X}$ denote by $C(\mathcal{X})$ the space of bounded continuous 
functions on $\mathcal{X}$. Let $\mathcal{M}(\K)$ be the vector space of finite signed Borel measures on $\K$
(and so $\M(\K)$ is its positive cone). Let
$\pi:\mathcal{M}(\K)\to\mathcal{M}(\Y)$ be defined by $(\pi\mu)(B)=\mu((\R^n\times B)\cap\K)$ for all $B\in\mathcal{B}(\Y)$, with adjoint mapping $\pi^*:C(\Y)\to C(\K)$ defined as
\[(\bx,\y)\,\mapsto\,(\pi^*h)(\bx,\y)\,:=\,h(\y),\qquad \forall h\in C(\Y).\]
Put (\ref{pb2}) in the framework of infinite-dimensional linear programs on vector spaces,
as described in e.g. \cite{anderson}. That is:
\[\rho=\inf_{\mu\in\mathcal{M}(\K)} \{\langle f,\mu\rangle\::\:\pi\mu=\varphi,\:\mu\geq0\},\]
with dual:
\[\tilde{\rho}=\sup_{h\in C(\Y)} \{\langle h,\varphi\rangle\::\:f-\pi^*h\geq 0\quad\mbox{on }\K\}.\]
One first proves that $\rho=\tilde{\rho}$ and then $\tilde{\rho}=\rho^*$. 

By \cite[Theor. 3.10]{anderson}, to get $\rho=\tilde{\rho}$, it suffices to prove that the set
$D:=\{(\pi\mu,\langle f,\mu\rangle):\mu\in\M(\K)\}$ is closed for the respective weak $\star$ topologies
$\sigma(\mathcal{M}(\Y)\times\R,C(\Y)\times\R)$ and $\sigma(\mathcal{M}(\K),C(\K))$.
Therefore consider a  converging sequence $\pi\mu_n\to a$ with $\mu_n\in\M(\K)$. The
sequence $(\mu_n)$ is uniformly bounded because 
\[\mu_n(\K)=(\pi\mu_n)(\Y)\,=\,\langle {\rm 1},\pi\mu_n\rangle\,\to\,\langle {\rm 1},a\rangle\,=\, a(\Y).\]
But by the Banach-Alaoglu Theorem (see e.g. \cite{ash}), the bounded closed sets of $\M(\K)$ are 
compact in the weak $\star$ topology. And so $\mu_{n_k}\to\mu$ for some $\mu\in\M(\K)$ and some subsequence $(n_k)$. Next, observe that for  $h\in C(\Y)$ arbitrary,
\[\langle h,\pi\mu_{n_k}\rangle\,=\,
\langle \pi^*h,\mu_{n_k}\rangle\,\to\,
\langle \pi^*h,\mu\rangle\,=\,\langle h,\pi\mu\rangle,\]
where we have used that $\pi^*h\in C(\K)$.
Hence combining the above with $\pi\mu_{n_k}\to a$, we obtain $\pi\mu=a$.
Similarly, $\langle f,\mu_{n_k}\rangle\to\langle f,\mu\rangle$ because $f\in C(\K)$.
Hence $D$ is closed and the desired result $\rho=\tilde{\rho}$ follows.

We next prove that $\tilde{\rho}=\rho^*$. Given $\epsilon>0$ fixed arbitrary, there is a function $h_\epsilon\in C(\Y)$ such that
$f-h_\epsilon\geq 0$ on $\K$ and $\int h_\epsilon d\varphi\geq\tilde{\rho}-\epsilon$.
By compactness of $\Y$ and the Stone-Weierstrass theorem,
 there is $p_\epsilon\in\R[\y]$ such that $\sup_{\y\in\Y}\vert h_\epsilon(\y)-p_\epsilon(\y)\vert\leq \epsilon$.
 Hence the polynomial $\tilde{p}_\epsilon:=p_\epsilon-\epsilon$ is feasible with value
 $\int_\Y\tilde{p}_\epsilon d\varphi\geq \tilde{\rho}-3\epsilon$, and as $\epsilon$ was arbitrary,
 the result $\tilde{\rho}=\rho^*$  follows.
\end{proof}

As next shown, optimal or nearly optimal solutions of $\P^*$ provide us with polynomial lower approximations of the optimal value function $\y\mapsto J(\y)$
that converges to $J(\cdot)$ in the $L_1(\varphi)$ norm. Moreover,
one may also obtain a piecewise polynomial approximation that converges to $J(\cdot)$ almost uniformly. (Recall that a sequence of measurable functions $(g_n)$ on a measure space 
$(\Y,\mathcal{B}(\Y),\varphi)$ converges to $g$ {\it almost uniformly}
if and only if for every $\epsilon>0$, there is a set $A\in\mathcal{B}(\Y)$ such that $\varphi(A)<\epsilon$ and
$g_n\to g$ uniformly on $\Y\setminus A$.)

\begin{cor}
\label{cor-dual}
Let both $\Y\subset\R^p$ and $\K$ in (\ref{set-xy}) 
be compact and assume that for every $\y\in\Y$,
the set $\K_\y$ is nonempty. Let $\P^*$ be as in (\ref{dual-lp}). 
If $(p_i)_{i\in\N}\subset\R[\y]$ is a maximizing sequence of (\ref{dual-lp}) 
then 
\begin{equation}
\label{cor-dual-1}
\displaystyle \int_\Y\,\vert\,J(\y)-p_i(\y)\,\vert\,d\varphi\,\to\,0\quad\mbox{as $i\to\infty$}.
\end{equation}
Moreover, define the functions $(\tilde{p}_i)$ as follows:
\[\tilde{p}_0:=p_0,\quad \y\mapsto\tilde{p}_i(\y)\,:=\,\max\,[\,\tilde{p}_{i-1}(\y),p_i(\y)\,],\quad i=1,2,\ldots\]
Then $\tilde{p}_i\to J(\cdot)$ almost uniformly.
\end{cor}
\begin{proof}
By Lemma \ref{lem-dual}, we already know that $\rho^*=\rho$ and so
\[\int_\Y p_i(\y)\,d\varphi(\y)\,\uparrow\,\rho^*\,=\,\rho\,=\,\int_\Y J(\y)\,d\varphi.\]
Next by feasibility of $p_i$ in (\ref{dual-lp}) 
\[f(\bx,\y)\geq p_i(\y)\quad\forall (\bx,\y)\in\K\:\Rightarrow\:
\inf_{\bx\in\K_\y}f(\bx,\y)=J(\y)\,\geq\,p_i(\y)\quad\forall\,\y\in\Y.\]
Hence (\ref{cor-dual-1}) follows from $p_i(\y)\leq J(\y)$ on $\Y$.

Next, with $\y\in\Y$ fixed, the sequence $(\tilde{p}_i(\y))_i$ is obviously monotone non decreasing and bounded above by $J(\y)$, hence with a limit $p^*(\y)\leq J(\y)$.
Therefore $\tilde{p}_i$ has the pointwise limit $\y\mapsto p^*(\y)\leq J(\y)$.
Also, by the Montone convegence theorem,  
$\int_\Y \tilde{p}_i(\y)d\varphi(\y)\to \int_\Y p^*(\y)d\varphi(\y)$. This latter fact
combined with (\ref{cor-dual-1}) and $p_i(\y)\leq \tilde{p}_i(\y)\leq J(\y)$ yields
\[0\,=\,\int_\Y (J(\y)-p^*(\y))\,d\varphi(\y),\]
which in turn implies that $p^*(\y)=J(\y)$ for almost all $\y\in\Y$.
Therefore $\tilde{p}_i(\y)\to J(\y)$ for almost all $\y\in\Y$.
And so, by Egoroff's Theorem \cite[Theor. 2.5.5]{ash}, $\tilde{p}_i\to J(\cdot)$ almost uniformly.
\end{proof}

\section{A hierarchy of semidefinite relaxations}

In general, solving the infinite-dimensional problem $\P$ and getting an optimal solution $\mu^*$ is impossible.
One possibility is to use numerical discretization schemes on a box containing $\K$; see for instance
\cite{ohl-lp}. But in the present context  of parametric optimization, 
if one selects finitely many {\it grid} points $(\bx,\y)\in\K$,
one is implicitly considering solving (or rather approximating) $\P_\y$ for finitely many points $\y$ in a grid of $\Y$, 
which we want to avoid.
To avoid this numerical discretization scheme we will use specific features of $\P$ when its data 
$f$ (resp. $\K$) is a polynomial (resp. a compact basic semi-algebraic set).

Therefore in this section we are now considering a {\it polynomial} parametric optimization problem, a special case of (\ref{pb1}) as we assume the following:

\begin{itemize}
\item $f\in\R[\bx,\y]$ and $h_j\in\R[\bx,\y]$, for every $j=1,\ldots,m$.
\item $\K$ is compact and $\Y\subset\R^p$ is a compact basic semi-algebraic set.
\end{itemize}
Hence the set $\K\subset\R^n\times\R^p$ in (\ref{set-xy}) is a compact basic semi-algebraic set.
We also assume that there is a probability measure $\varphi$ on $\Y$, absolutely continuous with respect to the Lebesgue measure, whose moments $\gamma=(\gamma_\beta)$, $\beta\in\N^p$, are available.
As already mentioned, if $\Y$ is a simple set (like e.g. a simplex or a box) then one may choose $\varphi$ to be the probability measure uniformly distributed on $\Y$, for which all moments can be computed easily. 
Sometimes, in the context of optimization with data uncertainty, the probability measure $\varphi$ is already specified
and in this case we assume that its moments $\gamma=(\gamma_\beta)$, $\beta\in\N^p$, are available.

\subsection{Notation and preliminaries}
Let $\N^n_i:=\{\alpha\in\N^n:\vert\alpha\vert\leq i\}$ with $\vert\alpha\vert=\sum_i\alpha_i$.
With a sequence $\z=(z_{\alpha\beta})$, $\alpha\in\N^n,\beta\in\N^p$, indexed in the canonical basis
$(\bx^\alpha\,\y^\beta)$  of $\R[\bx,\y]$, let 
$L_\z:\R[\bx,\y]\to\R$ be the linear mapping:
\[f\:(=\sum_{\alpha\beta}f_{\alpha\beta}(\bx,\y))\,\mapsto\: L_\z(f)\,:=\,
\sum_{\alpha\beta}f_{\alpha\beta}\,z_{\alpha\beta},\qquad f\in\R[\bx,\y].\]

\subsubsection*{Moment matrix}
The moment matrix $\M_i(\z)$ associated with a sequence $\z=(z_{\alpha\beta})$, has its rows and columns
indexed in the canonical basis $(\bx^\alpha\,\y^\beta)$, and with entries.
\[\M_i(\z)(\alpha,\beta),(\delta,\gamma))\,=\,L_\z(\bx^{\alpha}\y^\beta\,x^\delta\y^\gamma)\,=\,z_{(\alpha+\delta)(\beta+\gamma)},\]
for every $\alpha,\delta\in\N^n_i$ and every $\beta,\gamma\in\N^p_i$.

\subsubsection*{Localizing matrix}
Let $q$ be the polynomial $(\bx,\y)\mapsto q(\bx,\y)
:=\sum_{u,v}q_{uv}\bx^u\y^v$. The localizing matrix $\M_i(q\,\z)$ associated with 
$q\in\R[\bx,\y]$ and a sequence
$\z=(z_{\alpha\beta})$, has its rows and columns
indexed in the canonical basis $(\bx^\alpha\,\y^\beta)$, and with entries.
\begin{eqnarray*}
\M_i(q\,\z)(\alpha,\beta),(\delta,\gamma))&=&L_\z(q(\bx,\y)\bx^{\alpha}\y^\beta\,x^\delta\y^\gamma)\\
&=&\sum_{u\in\N^n,v\in\N^p}q_{uv}z_{(\alpha+\delta+u)(\beta+\gamma+v)},\end{eqnarray*}
for every $\alpha,\delta\in\N^n_i$ and every $\beta,\gamma\in\N^p_i$.

A sequence $\z=(z_{\alpha\beta})\subset\R$ has a {\it representing} finite Borel measure supported on $\K$
if there exists a finite Borel measure $\mu$ such that
\[z_{\alpha\beta}\,=\,\int_\K \bx^\alpha\,\y^\beta\,d\mu,\qquad\forall\,\alpha\in\N^n,\,\beta\in\N^p.\]
The next important result states a necssary and sufficient condition when $\K$
is compact and its defining polynomials $(h_k)\subset\R[\bx,\y]$ satisfy some condition.

\begin{assumption}
\label{assput}
Let $(h_j)_{j=1}^t\subset\R[\bx,\y]$ be a given family of polynomials.
There is some $N$ such that
the quadratic polynomial $(\bx,\y)\mapsto N-\Vert (\bx,\y)\Vert^2$ can be written
\[N-\Vert (\bx,\y)\Vert^2\,=\,\sigma_0+\sum_{k=1}^t\sigma_j\,h_j,\]
for some s.o.s. polynomials $(\sigma_j)_{j=1}^t\subset\Sigma[\bx,\y]$. 
\end{assumption}
\begin{thm}
\label{thput}
Let $\K:=\{(\bx,\y)\,:\,h_k(\bx,\y)\geq0,\:j=1,\ldots t\}$ and let $(h_k)_{k=1}^t$ satisfy 
Assumption \ref{assput}.
A sequence $\z=(z_{\alpha\beta})$
has a representing measure on $\K$ if and only if:
\[\M_i(\z)\succeq\,0\,;\quad\M_i(h_k\,\z)\,\succeq\,0,\quad k=0,\ldots,t.\]
\end{thm}
Theorem \ref{thput} is a direct consequence of Putinar's Positivstellensatz \cite{putinar} and \cite{markus}.
Of course, when Assumption \ref{assput} holds then $\K$ is compact. On the other hand, if $\K$ is compact and 
one knows a bound $N$ for $\Vert(\bx,\y)\Vert$ on $\K$ then its suffices to add the redundant 
quadratic constraint $h_{t+1}(\bx,\y)(:=N^2-\Vert(\bx,\y)\Vert^2)\geq0$ to the definition of $\K$,
and Assumption \ref{assput} holds.

\subsection{Semidefinite relaxations}

To compute (or at least, approximate)
the optimal value $\rho$ of problem $\P$ in (\ref{pb2}),
we now provide a hierarchy of semidefinite relaxations in the spirit of 
those defined in \cite{lassiopt}.

Let $\K\subset\R^n\times\R^p$ be as in (\ref{set-xy}), and let
$\Y\subset\R^p$ be the compact semi-algebraic set defined by:
\begin{equation}
\label{sety}
\Y\,:=\,\{\,\y\in\R^p\::\:h_k(\y)\,\geq\,0,\quad k=m+1,\ldots,t\}\end{equation}
for some polynomials $(h_k)_{k=m+1}^t\in\R[\y]$;
let $v_k:=\lceil ({\rm deg}\,h_k)/2\rceil]$ for every $k=1,\ldots,t$.
Next, let $\gamma=(\gamma_\beta)$ with
\[\gamma_\beta\,=\,\int_\Y\y^\beta\,d\varphi(\y),
\qquad\forall\,\beta\in\N^p,\]
be the moments of a probability  measure $\varphi$ on $\Y$,
absolutely continuous with respect to the Lebesgue measure, 
and let $i_0:=\max[\lceil ({\rm deg}\, f)/2\rceil,\max_kv_k]$.
For $i\geq i_0$, consider the following semidefinite relaxations:

\begin{equation}
\label{primal}
\begin{array}{rll}
\rho_i=&\displaystyle\inf_\z& L_\z(f)\\
&\mbox{s.t.}&\M_i(\z)\succeq0\\
&&\M_{i-v_j}(h_j\,\z)\succeq0,\quad j=1,\ldots,t\\
&&L_\z(\y^\beta)=\gamma_\beta,\quad \forall\,\beta\in\N^p_i.\end{array}
\end{equation}
\begin{thm}
\label{th2}
Let $\K,\Y$ be as (\ref{set-xy}) and (\ref{sety}) respectively, and let
$(h_k)_{k=1}^t$ satisfy Assumption \ref{assput}.
Assume that for every $\y\in\Y$ the set $\K_\y$ is nonempty, 
and for almost all $\y\in\Y$, $J(\y)$ is attained at a unique optimal solution 
$\bx^*(\y)$. Consider the semidefinite
relaxations (\ref{primal}). Then:

{\rm (a)} $\rho_i\uparrow\rho$ as $i\to\infty$.

{\rm (b)} Let $\z^i$ be a nearly optimal solution of (\ref{primal}), e.g. such that
$L_{\z^i}(f)\leq\rho_i+1/i$, and let $g\,:\Y\to\,\K_\y$ be the measurable mapping in Theorem \ref{th1}(c). 
Then 
\begin{equation}
\label{th2-2}
\displaystyle\lim_{i\to\infty}z^i_{\alpha\beta}\,=\,
\int_\Y \y^\beta \,g(\y)^\alpha\,d\varphi(\y),\qquad\forall\,\alpha\in\N^n,\,\beta\in\N^p.\end{equation}
In particular, for every $k=1,\ldots,n$, 
\begin{equation}
\label{th2-1}
\displaystyle\lim_{i\to\infty}z^i_{e(k)\beta}\,=\,
\int_\Y \y^\beta \,g_k(\y)\,d\varphi(\y),\qquad\forall\,\beta\in\N^p,\end{equation}
where $e(k)=(\delta_{j=k})_j\in\N^n$.
\end{thm}
The proof is postponed to Section \ref{proofs}.
\begin{rem}
{\rm Observe that if $\rho_i=+\infty$ for some index $i$ in the hierarchy
(and hence for all $i'\geq i$),
then the set $\K_\y$ is empty for
all $\y$ in some Borel set of $\Y$ with $\varphi(\Y)>0$. 
Conversely, one may prove that if $\K_\y$ is empty for
all $\y$ in some Borel set of $\Y$ with $\varphi(\Y)>0$, then necessarily 
$\rho_i=+\infty$ for all $i$ sufficiently large. In other words, the 
hierarchy of semidefinite relaxations (\ref{primal}) may also provide a certificate of
emptyness of $\K_y$ for some Borel set of $\Y$ with positive Lebesgue measure.}
\end{rem}
\subsection{The dual semidefinite relaxations}

The dual of the semidefinite relaxtion (\ref{primal}) reads:
\begin{equation}
\label{dual}
\begin{array}{rl}
\rho^*_i=&\displaystyle\sup_{p,(\sigma_i)} \int_\Y p\,d\varphi\\
\mbox{s.t.}&f-p=\sigma_0+\sum_{j=1}^t\sigma_j\,h_j\\
&\\
&p\in\R[\y];\:\sigma_j\subset\Sigma[\bx,\y],\quad j=1,\ldots,t\\
&{\rm deg}\,p\leq 2i,\,{\rm deg}\,\sigma_jh_j\leq 2i,\quad j=1,\ldots,t\\
\end{array}
\end{equation}
Observe that (\ref{dual}) is a strenghtening of (\ref{dual-lp}) as one restricts to polynomials $p\in\R[\y]$ of degree at most $2i$
and the nonnegativity of $f-p$ in (\ref{dual-lp}) is replaced with a stronger requirement in (\ref{dual}).
Therefore $\rho^*_i\leq\rho^*$ for every $i$.

\begin{thm}
\label{th22}
Let $\K,\Y$ be as (\ref{set-xy}) and (\ref{sety}) respectively, and let
$(h_k)_{k=1}^t$ satisfy Assumption \ref{assput}.
Assume that for every $\y\in\Y$ the set $\K_\y$ is nonempty,
and consider the semidefinite
relaxations (\ref{dual}). Then:

{\rm (a)} $\rho^*_i\uparrow\rho$ as $i\to\infty$.

{\rm (b)} Let $(p_i,(\sigma_j^i))$ be a nearly optimal solution of (\ref{dual}), e.g. such that
$\int_\Y p_id\varphi \geq\rho^*_i-1/i$. Then $p_i\leq J(\cdot)$ and
\begin{equation}
\label{th22-2}
\displaystyle\lim_{i\to\infty}\int_\Y(J(\y)-p_i(\y))\,d\varphi(\y)=0
\end{equation}
Moreover if one defines
\[\tilde{p}_0:=p_0,\quad \y\mapsto\tilde{p}_i(\y)\,:=\,\max\,[\,\tilde{p}_{i-1}(\y),p_i(\y)\,],\quad i=1,2,\ldots,\]
then $\tilde{p}_i\to J(\cdot)$ almost uniformly on $\Y$.
\end{thm}
\begin{proof}
Recall that by Lemma \ref{lem-dual}, $\rho=\rho^*$. Moreover 
let $(p_k)\subset\R[\y]$ be a maximizing sequence of (\ref{dual-lp}) as in Corollary
\ref{cor-dual} with value $s_k:=\int p_kd\varphi$, and let $p'_k:=p_k-1/k$ for every $k$ so that 
$f-p'_k>1/k$ on $\K$. By Theorem \ref{thput}, there exist s.o.s. polynomials $(\sigma_j^k)\subset\Sigma[\bx,\y]$ such that
$f-p'_k=\sigma_0^k+\sum_j\sigma_j^kh_j$. Letting $d_k$ be the maximum degree of $\sigma_0$ and $\sigma_jh_j$, $j=1,\ldots,t$,
it follows that $(s_k-1/k,(\sigma_j^k))$ is a feasible solution of (\ref{dual}) with $i:=d_k$. Hence 
$\rho^*\geq \rho^*_{d_k}\geq s_k-1/k$ and the result (a) follows because $s_k\to \rho^*$, and the sequence $\rho^*_i$ is monotone.
Then
(b) follows from Corollary \ref{cor-dual}.

\end{proof}
Hence Theorem \ref{th22} provides a lower polynomial approximation $p_i\in\R[\y]$ of the optimal value function $J(\cdot)$.
Its degree is bounded by $2i$, the order of the moments $(\gamma_\beta)$ of $\varphi$ taken into account
in the semidefinite relaxation (\ref{dual}). Moreover one may even define
a piecewise polynomial lower approximation $\tilde{p}_i$ that converges almost uniformly to $J(\cdot)$ on $\Y$.

\subsection*{Functionals of the optimal solutions}

Theorem \ref{th2} provides a mean of approximating any polynomial functional
on the optimal solutions of $\P_\y$, $\y\in\Y$. Indeed,

\begin{cor}
\label{coro2}
Let $\K,\Y$ be as (\ref{set-xy}) and (\ref{sety}) respectively, and let
$(h_k)_{k=1}^t$ satisfy Assumption \ref{assput}.
Assume that for every $\y\in\Y$ the set $\K_\y$ is nonempty, 
and for almost all $\y\in\Y$, $J(\y)$ is attained at a unique optimal solution 
$\bx^*(\y)\in\X^*_\y$. 
Let $h\in\R[\bx]$,
\[\bx\mapsto h(\bx)\,:=\,\sum_{\alpha\in\N^n}h_\alpha\,\bx^\alpha,\]
and let $\z^i$ be a nearly optimal solution of
the semidefinite relaxations (\ref{primal}).  

Then, for $i$ sufficiently large,
\[\int_\Y h(\bx^*(\y))\,d\varphi(\y)\,\approx\,\sum_{\alpha\in\N^n}h_\alpha\,z^i_{\alpha0}.\]
\end{cor}
\begin{proof}
The proof is an immediate consequence of Theorem \ref{th2} and Corollary \ref{coro1}.
\end{proof}

\subsection{Persistence for Boolean variables}

One interesting and potentially useful application is in Boolean optimization. Indeed suppose that 
for some subset $I\subseteq\{1,\ldots,n\}$, the variables $(x_i)$, $i\in I$, are boolean, that is,
the definition of $\K$ in (\ref{set-xy}) includes the constraints $x_i^2-x_i=0$, for every $i\in I$.

Then for instance, one might be interested to determine whether in an optimal solution $\bx^*(\y)$ of $\P_\y$,
and for some index $i\in I$, one has $x_i^*(\y)=1$ (or $x_i^*(\y)=0$)
for almost all values of the parameter $\y\in\Y$. 
In \cite{persis,natarajan} the probability that $x^*_k(\y)$ is $1$ is called the {\it persistency} 
of the boolean variable $x^*_k(\y)$

\begin{cor}
Let $\K,\Y$ be as in (\ref{set-xy}) and (\ref{sety}) respectively. Let $(h_k)_{k=1}^t$
satisfy (\ref{assput}).
Assume that for every $\y\in\Y$ the set $\K_\y$ is nonempty. 
Let $\z^i$ be a nearly optimal solution of
the semidefinite relaxations (\ref{primal}). Then for $k\in I$ fixed. 

{\rm (a)} $x^*_k(\y)=1$ for almost all $\y\in \Y$, only if $\displaystyle\lim_{i\to\infty}\,z^i_{e(k)0}=1$.

{\rm (b)} $x^*_k(\y)=0$ for almost all $\y\in \Y$, only if $\displaystyle\lim_{i\to\infty}\,z^i_{e(k)0}=0$.

Assume that for almost all $\y\in\Y$, $J(\y)$ is attained at a unique optimal solution $\bx^*(\y)\in\X^*_\y$.
Then ${\rm Prob}\,(x^*_k(\y)=1)=\displaystyle\lim_{i\to\infty}\,z^i_{e(k)0}$, and so:

{\rm (c)} $x^*_k(\y)=1$ for almost all $\y\in \Y$, if and only if $\displaystyle\lim_{i\to\infty}\,z^i_{e(k)0}=1$.

{\rm (d)} $x^*_k(\y)=0$ for almost all $\y\in \Y$, if and only if $\displaystyle\lim_{i\to\infty}\,z^i_{e(k)0}=0$.
\end{cor}
\begin{proof}
(a) The {\it only if part}. Let $\alpha:=e(k)\in\N^n$. From the proof of Theorem \ref{th2}, there is a subsequence $(i_l)_l\subset (i)_i$
such that 
\[\lim_{l\to\infty}z^{i_l}_{e(k)0}\,=\,\int_\K x_k\,d\mu^*,\]
where $\mu^*$ is an optimal solution of $\P$. Hence, by Theorem \ref{th1}(b), 
$\mu^*$ can be disintegrated into $\psi^*(d\bx\vert\y)d\varphi(\y)$ where 
$\psi^*(\cdot\vert\y)$ is a probability measure on $\X^*_\y$ for every $\y\in\Y$. Therefore, 
\begin{eqnarray*}
\lim_{l\to\infty}z^{i_l}_{e(k)0}&=&\int_\Y \left(\int_{\X^*_\y} x_k\psi^*(d\bx\,\vert\,\y)\right)\,d\varphi(y),\\
&=&\int_\Y \psi^*(\X^*_\y\,\vert\,\y)\,d\varphi(y)\quad\mbox{[because $x^*_k=1$ in $\X^*_\y$]}\\
&=&\int_\Y d\varphi(y)\,=\,1,
\end{eqnarray*}
and as the subsequence $(i_l)_l$ was arbitrary, the whole sequence $(z^i_{e(k)0})$ converges to $1$, the desired result.
The proof of (b) being exactly the same is omitted.

Next, if for every $\y\in \Y$, $J(\y)$ is attained at a singleton, by Theorem \ref{th2}(b),
\begin{eqnarray*}
\lim_{i\to\infty}z^i_{e(k)0}&=&\int _\Y x^*_k(\y)\,d\varphi(y)\,=\,\varphi(\{\y\,:\,x^*_k(\y)=1\})\\
&=&{\rm Prob}\,(x^*_k(\y)=1),\end{eqnarray*}
from which (c) and (d) follow.
\end{proof}

\subsection{Estimating the density $g(\y)$}

By Corollary \ref{coro2}, one may approximate any polynomial functional 
of the optimal solutions, like for instance the mean, variance, etc .. (with respect to the probability
measure $\varphi$). However, one may also wish to approximate (in some sense) the 
"curve" $\y\mapsto g_k(\y)$, that is, the surface described by the $k$-th coordinate 
$\bx^*_k(\y)$
of the optimal solution $\bx^*(\y)$ when $\y$ varies in $\Y$.

So let $g:\Y\to\R^n$ be the measurable mapping in Theorem \ref{th2} and 
suppose that one knows some lower bound vector $\a=(a_k)\in\R^n$, where:
\[a_k\,\leq\,\inf \:\{\:x_k\::\:(\bx,\y)\,\in\,\K\:\},\qquad k=1,\ldots,n.\]
Then for every $k=1,\ldots,n$, the measurable function $\hat{g}_k:\Y\to\R^n$ defined by
\begin{equation}
\label{hatgk}
\y\:\mapsto\quad \hat{g}_k(\y)\,:=\,g_k(\y)-a_k,\qquad \y\in\Y,\end{equation}
is nonnegative and integrable with respect to $\varphi$. 

Hence for every
$k=1,\ldots,n$, one may consider $d\lambda:=\hat{g}_kdx$ as a Borel measure
on $\Y$ with unknown density $\hat{g}_k$ with respect to $\varphi$, but with known moments
$\u=(u_\beta)$. Indeed, using (\ref{th2-1}),
\begin{eqnarray}
\nonumber
u_\beta\,:=\,\int_\Y \y^\beta\,d\lambda(\y)&=&-a_k\,\int_{\Y}\y^\beta\,d\varphi(\y)+\int_\Y \y^\beta\,g_k(\y)\,d\varphi(\y)\\
\label{momentu}
&=&-a_k\gamma_\beta\,+\,z_{e(k)\beta},\qquad\forall\beta\in\N^p,\end{eqnarray}
where for every $k=1,\ldots,n$, 
\[z_{e(k)\beta}\,=\,\displaystyle\lim_{i\to\infty}z^i_{e(k)\beta},\qquad \forall \beta\in\N^n,\]
with $\z^i$ being an optimal (or nearly optimal) solution of the semidefinite relaxation 
(\ref{primal}).

Hence we are now faced with a density estimation problem, that is:
Given the sequence of moments
$\gamma_\beta=\int_\Y \y^\beta g_k(\y)d\varphi$, $\beta\in\N^p$, of
the unknown nonnegative measurable function $g_k$ on $\Y$,  "estimate" $g_k$. 
One possibility is the so-called {\it maximum entropy} approach, briefly described in the next section.

\subsection*{Maximum-entropy estimation}

We briefly describe the {\it maximum entropy} estimation technique 
in the univariate case. The multivariate case generalizes easily.
Let $g\in L_1([0,1])$\footnote{$L_1([0,1])$ denote the Banach space of integrable functions 
on the interval $[0,1]$ of the real line, equipped with the norm $\Vert g\Vert_1 =\int_0^1\vert b(\bx)\vert\,d\bx$.} be a nonnegative function only known
via the first  $2d+1$ moments $\u=(u_j)_{j=0}^{2d}$ of its associated
measure $d\varphi=gd\bx$ on $[0,1]$. (In the context of previous section,
the function $g$ to estimate is $\y\mapsto g_k(\y)$ in (\ref{hatgk}) 
from the sequence $\u$ in (\ref{momentu}) of its (multivariate) moments.)

From that partial knowledge one wishes (a) 
to provide an estimate $h_d$ of $g$ such that
the first $2d+1$ moments of the measure $h_dd\bx$ match those of $gd\bx$, and (b) analyze 
the asymptotic behavior of $h_d$ when $d\to\infty$. This problem has important applications in various areas of physics, engineering, and signal processing in particular.

An elegant methodology is to search for $h_d$ in a (finitely) parametrized family $\{h_d(\lambda,x)\}$ of functions, and optimize over the unknown parameters $\lambda$ via a suitable criterion. For instance, 
one may wish to select an estimate $h_d$ that maximizes some appropriate {\it entropy}. 
Several choices of entropy functional are possible as long as one obtains a convex optimization problem in the finitely many coefficients $\lambda_i$'s. For more details the interested reader is referred to e.g.
Borwein and Lewis \cite{borwein1,borwein2} and the many references therein.

We here choose the Boltzmann-Shannon entropy $\H:L_1([0,1])\to  \R\cup\{-\infty\}$:
\begin{equation}
\label{shannon}
h\mapsto\,\H[h]\,:=\,-\int_0^1 h(x)\,\ln{h(x)}\,dx,
\end{equation}
a strictly concave functional. Therefore, the problem reduces to:
\begin{equation}
\label{chap14-3-pb1}
\sup_h\: \left\{\:\H[h]\::\:\int_0^1 x^j\,h(x)\,dx\,=\,u_j,\quad j=0,\ldots,2d\:\right\}.
\end{equation}
The structure of this infinite-dimensional convex optimization problem
permits to search for an optimal solution $h_d^*$ of the form:
\begin{equation}
\label{chap14-3-g*}
x\mapsto h_d^*(x)\,=\,\exp{\sum_{j=0}^{2d}\lambda_j^*\,x^j},
\end{equation}
and so $\lambda^*$ is an optimal solution of the finite-dimensional unconstrained convex problem
\[\theta(\u):=\sup_\lambda\:\langle\u,\lambda\rangle-
\int_0^1\exp\left(\sum_{j=0}^{2d}\lambda_jx^j\right)\,dx.\]
Notice that the above function $\theta$ is just the Legendre-Fenchel transform of the 
convex function $\lambda\mapsto \int_0^1\exp{\sum_{j=0}^{2d}\lambda_jx^j}\,dx$.

An optimal solution can be calculated by applying first-order
methods, in which case the gradient $\nabla v_d$ of the function
\[\lambda\mapsto v_d(\lambda):=\langle\u,\lambda\rangle-\int_0^1\exp\left(\sum_{j=0}^{2d}\lambda_jx^j\right)\,dx,\]
is provided by:
\[\frac{\partial v_d(\lambda)}{\partial \lambda_k}\,=\,u_k-\int_0^1x^k\,\exp\left(\sum_{j=0}^{2d}\lambda_jx^j\right)\,d\bx,\qquad k=0,\ldots,2d+1.\]
If one applies second-order methods, e.g. Newton's method, then computing the Hessian $\nabla^2v_d$ 
at current iterate $\lambda$, reduces to computing
\[\frac{\partial^2 v_d(\lambda)}{\partial \lambda_k\partial \lambda_j}
\,=\,-\int_0^1x^{k+j}\,\exp\left(\sum_{j=0}^{2d}\lambda_jx^j\right)\,dx,\qquad k,j=0,\ldots,2d+1.\]
In such simple cases like a box $[a,b]$ (or $[a,b]^n$ in the multivariate case)
such quantities can be approximated quite accurately via cubature formula
as described in e.g. \cite{gautschi}. In particular, several cubature formula behave very well
for exponentials of polynomials as shown in e.g. Bender et al. \cite{BMP}. An alternative 
with no cubature formula is also proposed in \cite{lascdc}.

One has the following convergence result which follows directly from
\cite[Theor. 1.7 and p. 259]{borwein1}.
\begin{prop}
Let $0\leq g\in L_1([0,1])$ and for every $d\in\N$, let $h^*_d$ in (\ref{chap14-3-g*}) be an optimal solution of
(\ref{chap14-3-pb1}). Then, as $d\to\infty$,
\[\int_0^1\psi(\y)\,(h^*_d(\y)-g(\y))\,dx\,\to\,0,\]
for every bounded measurable function $\psi:[0,1]\to\R$ which is continuous almost everywhere
\end{prop}
Hence, the max-entropy estimate we obtain is not a pointwise estimate of $g$, and so,
at some points of $[0,1]$ the max-entropy density $h^*_d$ and the density $g$
to estimate may differ significantly. However, for sufficiently large $d$, both curves of
$h^*_d$ and $g$ are close to each other. In our context, recall that 
$g$ is for instance $y\mapsto x^*_k(y)$, and so in general, for fixed $y$,
$h^*_d(y)$ is close to $x^*_k(y)$ and might be chosen for the $k$-coordinate 
of an initial point $\bx$, input of a local minimization algorithm to find the global minimizer $\bx^*(y)$.

\subsection{Illustrative examples}
In this section we provide some simple illustrative examples.
To show the potential of the approach we have voluntarily chosen very simple examples for which one knows the solutions exactly so as to compare the results we obtain with the exact 
optimal value and optimal solutions. The semidefinite relaxations (\ref{primal}) were implemented 
by using the software package Gloptipoly \cite{gloptipoly}.
The max-entropy estimate $h^*_d$ of $g_k$ was computed by 
using Newton's method, where at each iterate $(\lambda^{(k)},h_d(\lambda^{(k)}))$:
\[\lambda^{(k+1)}\,=\,\lambda^{(k)}-(\nabla^2v_d(\lambda^{(k)}))^{-1}\nabla v_d(\lambda^{(k)}).\]
\begin{ex}
\label{example1}
{\rm For illustration purpose, consider the toy example where $\Y:=[0,1]$,
\[\K:=\{(x,y)\::\:1-x^2+y^2\geq0;\,x,y\in\Y\}\subset\R^2,\quad (x,y)\mapsto f(x,y):=-x^2y.\]
Hence for each value of the parameter $y\in\Y$, the unique optimal solution is 
$x^*(y):=\sqrt{1-y^2}$. And so 
in Theorem \ref{th2}(b), $y\mapsto g(y)=\sqrt{1-y^2}$.

Let $\varphi$ be the probability measure uniformly distributed on $[0,1]$.
Therefore, 
\[\rho=\,\int_0^1J(y)\,d\varphi(y)\,=\,-\int_0^1y(1-y^2)\,dy\,=\,-1/4.\]
Solving (\ref{primal}) with $i:=3$, that is, with moments up to order $6$, one obtains the optimal value $-0.250146$. Solving (\ref{primal}) with $i:=4$,
one obtains the optimal value $-0.25001786$ and the moment sequence
\[\z\,=\,(1, 0.7812,0.5,0.6604,0.3334,0.3333,0.5813,0.25,0.1964,0.25,0.5244,0.2,0.1333,\]
\[0.1334,0.2,0.4810,0.1667, 0.0981,0.0833,0.0983,0.1667)\]
Observe that
\[z_{1k}-\int_0^1y^k\sqrt{1-y^2}\,dy\approx O(10^{-6}),\qquad k=0,\ldots 4,\]
\[z_{1k}-\int_0^1y^k\sqrt{1-y^2}\,dy\approx O(10^{-5}),\qquad k=5,6,7.\]
Using a max-entropy approach to approximate the density $y\mapsto g(y)$ on $[0,1]$, 
with the first $5$ moments $z_{1k}$, $k=0,\ldots,4$, we find that the 
optimal function $h^*_4$ in (\ref{chap14-3-g*}) 
is obtained with
\[\lambda^*=(-0.1564, 2.5316, -12.2194, 20.3835, -12.1867).\]
Both curves of $g$ and $h^*_4$ are displayed in Figure \ref{gversush}. 
Observe that with only $5$ moments, the max-entropy solution $h^*_4$
approximates $g$ relatively well, even if it differs significantly at some points.
Indeed, the shape of $h^*_4$ resembles very much that of $g$.
\begin{figure}[ht]
\begin{center}
\resizebox{0.8\textwidth}{!}
{\includegraphics{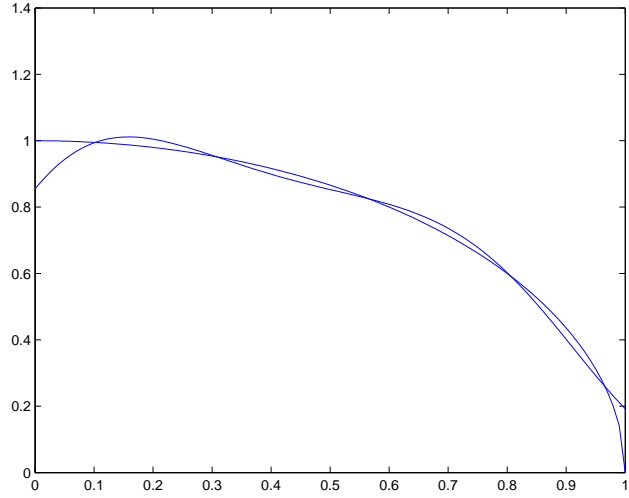}}
\caption{Example \ref{example1}: $g(y)=\sqrt{1-y^2}$ versus $h^*_4(y)$}
\label{gversush}
\end{center}
\end{figure}
 
Finally, from an optimal solution of (\ref{dual})
one obtains for $p\in\R[y]$, the degree-$8$ univariate polynomial
\begin{eqnarray*}
y&\mapsto&p(y)= -0.0004   -0.9909y   -0.0876y^2    +1.4364y^3   -1.2481y^4\\
&&+ 2.1261y^5   -2.1309y^6    +1.1593y^7   -0.2641y^8\end{eqnarray*}
and Figure \ref{diff0} displays the curve $y \mapsto J(y)-p(y)$ on $[0,1]$. One observes that
 $J\geq p$ and the maximum difference is about $3.10^{-4}$ close to $0$ and much less 
 for $y\geq 0.1$, a good precision with only $8$ moments.
 
\begin{figure}[ht]
\begin{center}
\resizebox{0.8\textwidth}{!}
{\includegraphics{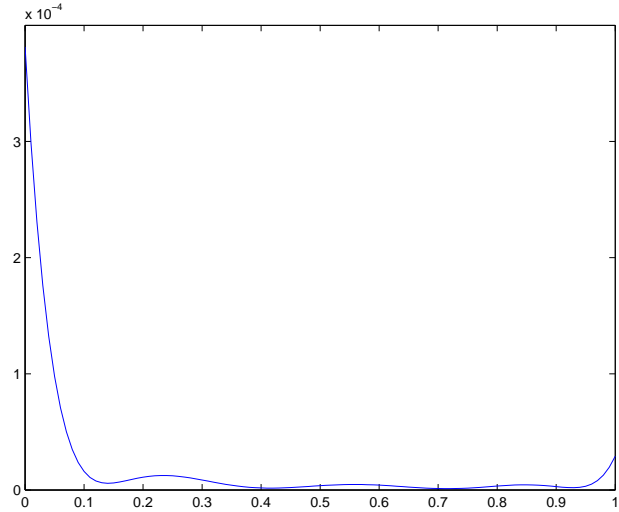}}
\caption{Example \ref{example1}: $J(y)-p(y)$ on $[0,1]$}
\label{diff0}
\end{center}
\end{figure}
}
\end{ex}
\begin{ex}
\label{example2}
{\rm Again with $\Y:=[0,1]$, let
\[\K:=\{(\bx,y)\::\:1-x_1^2-x_2^2\geq0\}\subset\R^2,\quad (\bx,y)\mapsto 
f(\bx,y):=y x_1+(1-y)x_2.\]
For each value of the parameter $y\in\Y$, the unique optimal solution $\bx^*\in\K$ satisfies
\[(x_1^*(y))^2+(x_2^*(y))^2=1;\quad (x^*_1(y))^2=\frac{y^2}{y^2+(1-y)^2},\quad (x^*_2(y))^2=\frac{(1-y)^2}{y^2+(1-y)^2},\]
with optimal value 
\[J(y)\,=\,-\frac{y^2}{\sqrt{y^2+(1-y)^2}}-\frac{(1-y)^2}{\sqrt{y^2+(1-y)^2}}\,=\,-\sqrt{y^2+(1-y)^2}.\]
So in Theorem \ref{th2}(b), 
\[y\mapsto g_1(y)=\frac{-y}{\sqrt{y^2+(1-y)^2}},\quad y\mapsto g_2(y)=\frac{y-1}{\sqrt{y^2+(1-y)^2}},\]
and with $\varphi$ being the probability measure uniformly distributed on $[0,1]$,
\[\rho=\,\int_0^1J(y)\,d\varphi(y)\,=\,-\int_0^1\sqrt{y^2+(1-y)^2}\,dy\,\approx\,-0.81162\]
Solving (\ref{primal}) with $i:=3$, that is, with moments up to order $6$, one obtains 
$\rho_3\approx-0.8117$ with $\rho_3-\rho\approx O(10^{-5})$.
Solving (\ref{primal}) with $i:=4$, one obtains 
$\rho_4\approx-0.81162$ with $\rho_4-\rho\approx O(10^{-6})$, and the moment sequence $(z_{k10})$, $k=0,1,2,3,4$:
\[z_{k10}\,=\,(-0.6232,\:   -0.4058,\:    -0.2971,\:    -0.2328,\:   -0.1907),\]
and 
\[z_{k10}-\int_0^1 y^k\,g_1(y)\,dy\,\approx\,O(10^{-5}),\quad k=0,\ldots, 4.\]

Using a max-entropy approach to approximate the density $y\mapsto -g_1(y)$ on $[0,1]$, 
with the first $5$ moments $z_{1k}$, $k=0,\ldots,4$, we find that the 
optimal function $h^*_4$ in (\ref{chap14-3-g*}) 
is obtained with
\[\lambda^*=(-3.61284,  15.66153266  -29.430901  27.326347   -9.9884452).\]
and we find that
\[z_{k10}\,+\,\int_0^1 y^k\,h^*_4(y)\,dy\,\approx\,O(10^{-11}),\quad k=0,\ldots,4.\]
In Figure \ref{gversush4} are displayed the two functions $-g_1$ and $h^*_4$, and one observes a very good
concordance.
\begin{figure}[ht]
\begin{center}
\resizebox{0.8\textwidth}{!}
{\includegraphics{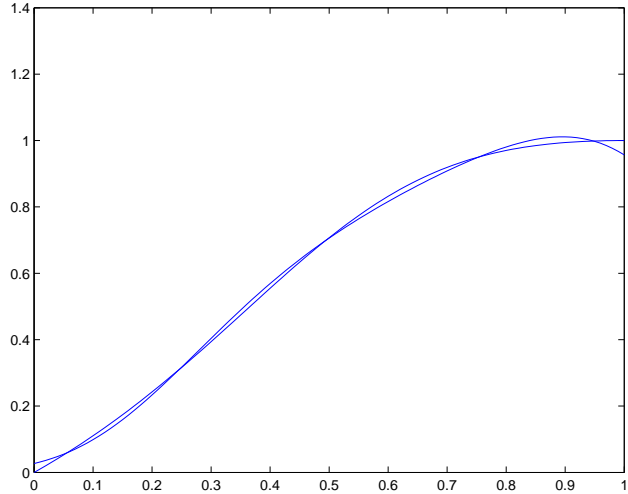}}
\caption{Example \ref{example2}: $h^*_4(y)$ versus $-g_1(y)=y/\sqrt{y^2+(1-y)^2}$}
\label{gversush4}
\end{center}
\end{figure}

Finally, from an optimal solution of (\ref{dual})
one obtains for $p\in\R[y]$, the degree-$8$ univariate polynomial
\begin{eqnarray*}
 x&\mapsto &p(y):=-1.0000   +0.9983y  -0.4537y^2  -0.9941y^3   +2.2488y^4   -7.6739y^5 \\
 &&  +11.8448y^6   -7.9606y^7   +1.9903y^8 \end{eqnarray*}
 and Figure \ref{diff} displays the curve $y \mapsto J(y)-p(y)$ on $[0,1]$. One observes that
 $J\geq p$ and the maximum difference is about $10^{-4}$, a good precision with only $8$ moments.
 
\begin{figure}[ht]
\begin{center}
\resizebox{0.8\textwidth}{!}
{\includegraphics{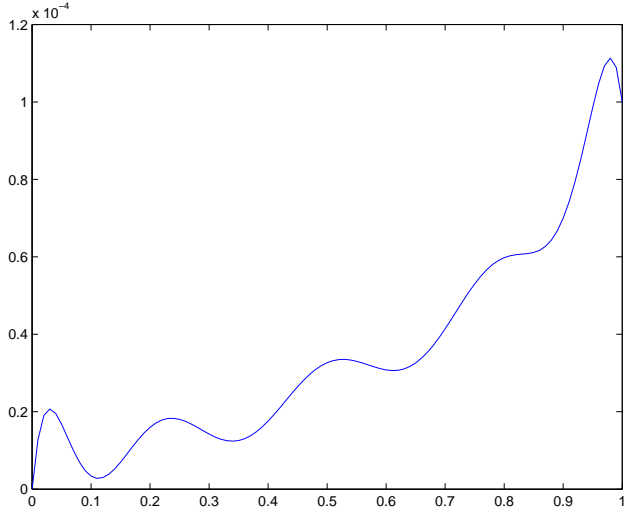}}
\caption{Example \ref{example2}: $J(y)-p(y)$ on $[0,1]$}
\label{diff}
\end{center}
\end{figure}
}
\end{ex}
\begin{ex}
\label{example33}
{\rm In this example one has $\Y=[0,1]$, $(\bx,y)\mapsto f(\bx,y):=y x_1+(1-y)x_2$, and
\[\K:=\{(\bx,y)\,:\: y x_1^2+x_2^2-y <=0;\,x_1^2+y x^2-y <=0\}.\]
That is, for each $y\in\Y$ the set $\K_y$ is the intersection of two ellipsoids.
It is easy to chack that $1+x^*_i(y)\geq0$ for all $y\in\Y$, $i:=1,1$.
With $i=4$ the max-entropy estimate $y\mapsto h^*_4(y)$ for $1+x^*_1(y)$ is obtained with
\[  \lambda^*=(-0.2894,    1.7192 , -19.8381,   36.8285,  -18.4828),\]
whereas the max-entropy estimate $y\mapsto h^*_4(y)$ for  $1+x^*_2(y)$ is obtained with
\[\lambda^*=(-0.1018,   -3.0928,    4.4068,    1.7096,  -7.5782).\]
Figure \ref{diff333} displays the curves of $x^*_1(y)$ and $x^*_2y)$, as well as the constraint $h_1(\bx^*(y),y)$. Observe that 
$h_1(\bx^*(y),y)\approx 0$ on $[0,1]$ which means that
for almost all $y\in [0,1]$, at an optimal solution $\bx^*(y)$,
the constraint $h_1\leq 0$ is saturated. Figure
\ref{diff334} displays the curves of $h_1(\bx^*(y),y)$
and $h_2(\bx^*(y),y)$.
\begin{figure}[ht]
\begin{center}
\resizebox{0.8\textwidth}{!}
{\includegraphics{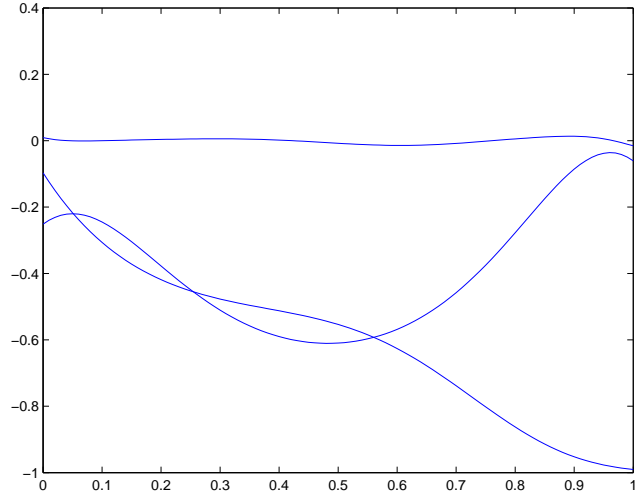}}
\caption{Example \ref{example33}: $x^*_1(y)$, $x^*_2(y)$ and $h_1(\bx^*(y),y)$ on $[0,1]$}
\label{diff333}
\end{center}
\end{figure}

\begin{figure}[ht]
\begin{center}
\resizebox{0.8\textwidth}{!}
{\includegraphics{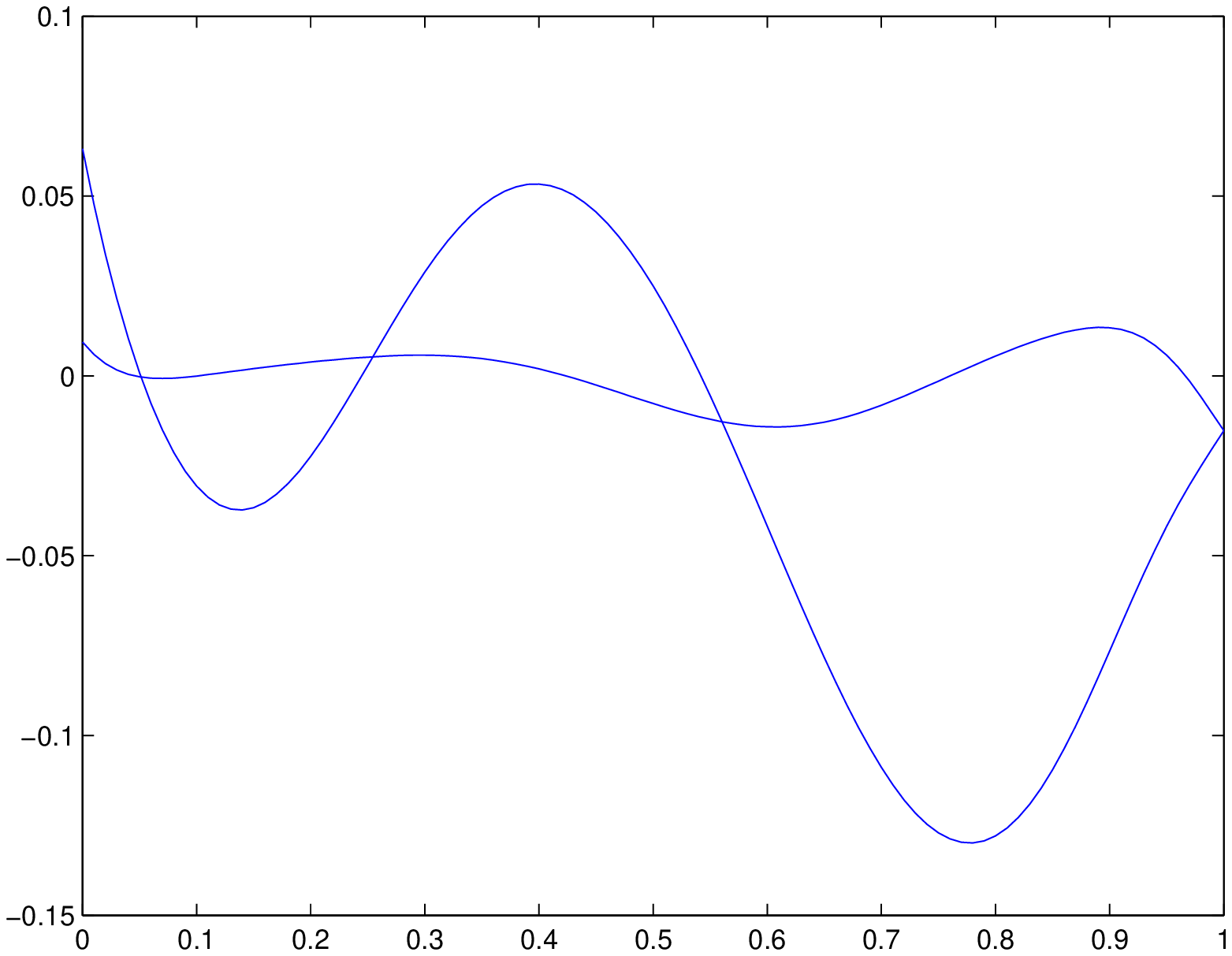}}
\caption{Example \ref{example33}: $h_1(\bx^*(y),y)$ and $h_2(\bx^*(y),y)$ on $[0,1]$}
\label{diff334}
\end{center}
\end{figure}
}

\end{ex}

\begin{ex}
\label{ex444}
{\rm This time $\Y=[0,1]$, $(\bx,y)\mapsto f(\bx,y):=(1-2y)(x_1+x_2)$, and
\[\K:=\{(\bx,y)\,:\: y x_1^2+x_2^2-y =0;\,x_1^2+y x^2-y =0\}.\]
That is, for each $y\in\Y$ the set $\K_y$ is the intersection of two ellipses, and
\[\bx=\left(\pm \sqrt{\frac{y}{1+y}},\pm \sqrt{\frac{y}{1+y}}\right);\quad J(y)=
-2\vert 1-2y\vert\,\sqrt{\frac{y}{1+y}}.\]
 With $i=4$ the max-entropy estimate $y\mapsto h^*_4(y)$ for $1+x^*_1(y)$ is obtained with
\[\lambda^*=(0.3071151, -12.51867,  43.215907, -46.985733,  16.395944).\]
In Figure \ref{diff444-cost} are displayed the curves $y\mapsto -p(y)$ and $y\mapsto -J(y)$,
whereas in Figure \ref{diff444-j-p} is displayed the curve $y\mapsto p(y)-J(y)$. One may
see that $p$ is a good lower approximation of $J$ even with only 8 moments.

\begin{figure}[ht]
\begin{center}
\resizebox{0.8\textwidth}{!}
{\includegraphics{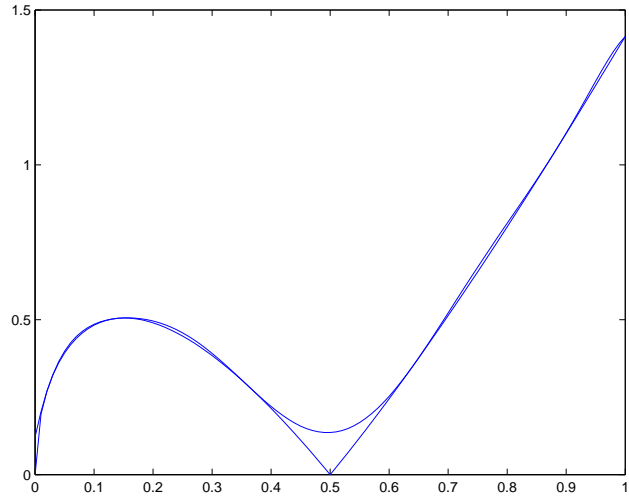}}
\caption{Example \ref{ex444}: $-p(y)$ and $-J(y)$ on $[0,1]$}
\label{diff444-cost}
\end{center}
\end{figure}
\begin{figure}[ht]
\begin{center}
\resizebox{0.8\textwidth}{!}
{\includegraphics{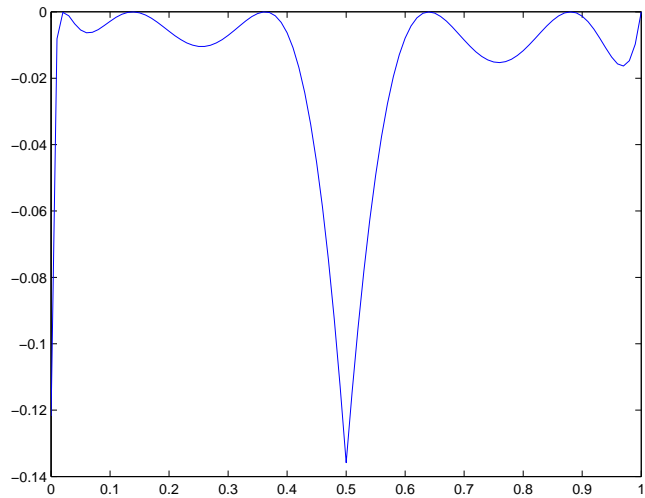}}
\caption{Example \ref{ex444}: the curve $p(y)-J(y)$ on $[0,1]$}
\label{diff444-j-p}
\end{center}
\end{figure}

On the other hand, in Figure \ref{diff444} is displayed $h^*_4(y)$ versus $x^*_1(y)$ where the latter is 
$-\sqrt{y/(1+y)}$ on $[0,1/2]$ and $\sqrt{y/(1+y)}$ on $[1/2,1]$. Here we see that
the discontinuity of $x^*_1(y)$ is difficult to approximate "pointwise" with few moments,
and despite a very good precision on the five first moments. Indeed:
\[\int_0^1 y^k\,(h^*_4(y)-1)\,dx-\int_0^1 y^k\,x^*_1(y)\,dx\,=\,O(10^{-14}),\quad k=0,\ldots,4.\]

\begin{figure}[ht]
\begin{center}
\resizebox{0.8\textwidth}{!}
{\includegraphics{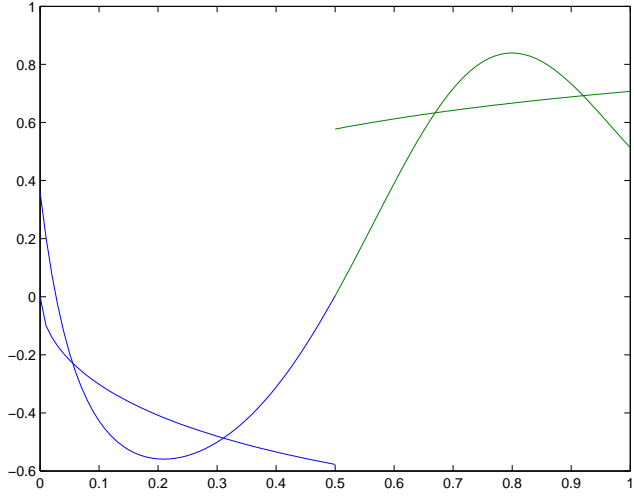}}
\caption{Example \ref{ex444}: $h^*_4(y)-1$ and $x^*_1(y)$ on $[0,1]$}
\label{diff444}
\end{center}
\end{figure}
}

\end{ex}
We end up this section with the case where the density $g_k$ to estimate is a step function
which would be the case in an optimization problem $\P_y$ with boolean variables
(e.g. the variable $x_k$ takes values in $\{0,1\}$).
\begin{ex}
\label{example3}
{\rm 
Assume that with a single parameter $y\in [0,1]$, 
the density $g_k$ to estimate is the step function. 
\[y\mapsto g_k(y):=\left\{\begin{array}{ll}
1&\mbox{if }y\in [0,1/3]\cup [2/3,1]\\
0&\mbox{otherwise.}\end{array}\right.\]
The max-entropy estimate $h^*_4$ in (\ref{chap14-3-g*}) with $5$ moments is
obtained with
\[\lambda^*=[-0.65473672   19.170724  -115.39354   192.4493171655  -96.226948865],\]
and we have 
\[\int_0^1y^kh^*_4(y)\,dy-\int_0^1y^k\,dg_k(y)\,\approx\,O(10^{-8}),\quad k=0,\ldots,4.\]
In particular, the persistency $\int_0^1 g_k(y)dy=2/3$ of the variable $x^*_k(y)$,
is very well approximated (up to $10^{-8}$ precision) by $\int h^*_4(y)dy$, with only $5$ moments.

\begin{figure}[ht]
\begin{center}
\resizebox{0.9\textwidth}{!}
{\includegraphics{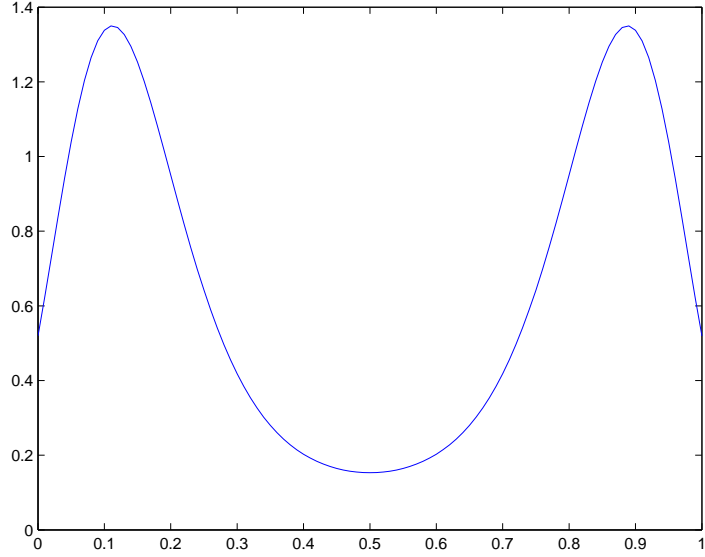}}
\caption{Example \ref{example3}: $g_k(y)={\rm 1}_{[0,1/3]\cup [2/3,1]}$ versus $h^*_4(y)$}
\label{gversush2}
\end{center}
\end{figure}
Of course, in this case and with only 5 moments, the density $h^*_4$ is not a good pointwise
approximation of the step function $g_k$; however its "shape" reveals the two steps of value $1$ separated by a step 
of value $0$. A better pointwise approximation would require more moments.
}
\end{ex}

\newpage

\section{Appendix}
\label{proofs}
\noindent
{\bf Proof of Theorem \ref{th2}.}~
\vspace{0.2cm}

 We already know that  $\rho_i\leq\rho$ for all $i\geq i _0$.
We also need to prove that   $\rho_i>-\infty$ for sufficiently large $i$. 
  Let  $Q\subset\R[\bx,\y]$ be the {\it quadratic module} generated by 
the polynomials $\{h_j\}\subset\R[\bx,\y]$ that define $\K$, i.e.,
\[Q\,:=\,\{\:\sigma\in\R[\bx,\y]\::\quad \sigma\,=\,\sigma_0+\sum_{j=1}^t\sigma_j\,h_j\quad \mbox{with }\{\sigma_j\}_{j=0}^t\subset\Sigma[\bx,\y]\}.\]
In addition, let  $Q(l)\subset Q$ be the set of elements $\sigma\in Q$ which have a representation $\sigma_0+\sum_{j=0}^t\sigma_j\,h_j$ for some s.o.s. family
 $\{\sigma_j\}\subset\Sigma^2$ with 
 ${\rm deg}\,\sigma_0\leq 2l$ and ${\rm deg}\,\sigma_jh_j\leq 2l$ for all $j=1,\ldots,t$.
 
 Let $i\in\N$ be fixed. As $\K$ is compact, there exists $N$ such that
 $N\pm \bx^\alpha\y^\beta >0$ on $\K$, for all $\alpha\in\N^n$ and $\beta\in\N^p$, with $\vert\alpha+\beta\vert\leq 2i$.
 Therefore, under Assumption \ref{assput}(ii), the polynomial  $N\pm \bx^\alpha\y^\beta$ belongs to $Q$; see Putinar \cite{putinar}.  But there is even some $l(i)$ such that
$N\pm \bx^\alpha\y^\beta\in Q(l(i))$ for every $\vert\alpha+\beta\vert\leq 2i$.
 Of course we also have $N\pm \bx^\alpha\y^\beta\in Q(l)$ for every $\vert\alpha+\beta\vert\leq 2i$, whenever  $l\geq l(i)$. Therefore, let us take $l(i)\geq i_0$.
  For every feasible solution $\z$ of $\Q_{l(i)}$ one has
  \[\vert z_{\alpha\beta}\vert\,=\,\vert\: L_{\z}(\bx^\alpha\y^\beta)\:\vert\leq N,\qquad \forall\,\vert\alpha+\beta\vert \leq 2i.\]
This follows from $z_0=1$, $\M_{l(i)}(\z)\succeq0$ and $\M_{l(i)-v_j}(h_j\,\z)\succeq0$, which implies
  \[Nz_0\pm z_{\alpha\beta}\,=\,L_{\z}(N\pm \bx^\alpha\y^\beta)\,=\,L_{\z}(\sigma_0)+\sum_{j=1}^tL_{\z}(\sigma_j\,h_j)\geq0\]
for some $\{\sigma_j\}\subset\Sigma[\bx,\y]$ with ${\rm deg}\,\sigma_j\,h_j\leq 2l(i)$. 
In particular, $L_\z(f)\,\geq\,-N\sum_{\alpha,\beta} \vert f_{\alpha\beta}\vert$, which proves that $\rho_{l(i)}>-\infty$, and so $\rho_i>-\infty$ for 
all sufficiently large $i$.\\

From what precedes, and with $k\in\N$ arbitrary, let $l(k)\geq k$ and $N_k$ be such that 
\begin{equation}
\label{newbound}
N_k\pm \bx^\alpha\y^\beta\in Q(l(k))\qquad \forall\,\alpha\,\in\N^n,\,\beta\in\N^p\:\mbox{ with }\vert\alpha+\beta\vert \leq 2k.
\end{equation}
Let $i\geq l(i_0)$, and let $\z^{i}$ be a nearly optimal solution of $(\ref{primal})$ with value
\begin{equation}
\label{value}
\rho_i\,\leq \,L_{\z^{i}}(f)\,\leq\,\rho_i+\frac{1}{i}\quad\left(\leq\rho+\frac{1}{i}\right).
\end{equation}
Fix $k\in\N$. Notice that from (\ref{newbound}), for every $i\geq l(k)$, one has 
\[\vert\,L_{\z^{i}}(\bx^\alpha\y^\beta)\,\vert\:\leq\: N_kz_0=N_k,\quad \forall\,\alpha\in\N^n,\beta\in\N^p \:\mbox{ with }
\vert\alpha+\beta\vert\,\leq\,2k.\]
Therefore, for all  $i\geq l(i_0)$,
\begin{equation}
\label{aux1}
\vert z^{i}_{\alpha\beta}\vert\:=\:\vert\,L_{\z^{i}}(\bx^\alpha\y^\beta)\,\vert\:\leq\: N'_k,\quad \forall\,\alpha\in\N^n,\beta\in\N^p \:\mbox{ with }\vert\alpha+\beta\vert\,\leq\,2k,
\end{equation}
where $N'_k=\max[N_k,V_k],$ with 
\[V_k:=\max_{\alpha,\beta, i}\:\{\:\vert z^{i}_{\alpha\beta}\vert\::\quad
\vert\alpha+\beta\vert\leq 2k\,;\quad l(i_0)\leq i\leq l(k)\:\}.\]

Complete each vector $\z^{i}$ with zeros to make it an infinite bounded sequence in $l_\infty$, 
indexed in the canonical basis $(\bx^\alpha\y^\beta)$ of $\R[\bx,\y]$. In view of (\ref{aux1}), 
\begin{equation}
\label{newbound2}
\vert z^{i}_{\alpha\beta}\vert\,\leq N'_k\qquad \forall\,\alpha\in\N^n,\beta\in\N^p\:\mbox{ with}\quad 2k-1\leq \vert\alpha+\beta\vert \leq 2k,
\end{equation}
and for all $k=1,2,\ldots$.

Hence, let $\widehat{\z}^{i}\in l_\infty$ be the new sequence defined by
\[\widehat{z}^{i}_{\alpha\beta}\,:=\, \frac{z^{i}_{\alpha\beta}}{N'_k},
\quad \forall\,\alpha\in\N^n,\beta\in\N^p\:\mbox{ with}\quad 2k-1\leq \vert\alpha+\beta\vert \leq 2k,\quad \forall\,k=1,2,\ldots,\]
and in $l_\infty$, consider the sequence $\{\widehat{\z}^{i}\}_i$, as $i\to\infty$.

Obviously, the sequence $\{\widehat{\z}^{i}\}_i$ is in the unit ball 
$B_1$ of $l_\infty$, and so,
by the Banach-Alaoglu theorem (see e.g. Ash \cite{ash}), there exists $\widehat{\z}\in B_1$, and a subsequence $\{i_l\}$, such that $\widehat{\z}^{i_l}\to \widehat{\z}$ 
as $l\to\infty$, for the weak $\star$ topology 
$\sigma(l_\infty,l_1)$ of $l_\infty$.
In particular, pointwise convergence holds, that is,
\[\lim_{l\to\infty}\,\widehat{z}^{i_l}_{\alpha\beta}\,\to\,\widehat{z}_{\alpha\beta}\qquad \forall\,\alpha\in\N^{n},\beta\in\N^p.\]
Next, define 
\[z_{\alpha\beta}\,:=\, \widehat{z}_{\alpha\beta}\times N'_k
\quad \forall\,\alpha\in\N^n,\beta\in\N^p\:\mbox{ with}\quad 2k-1\leq \vert\alpha+\beta\vert \leq 2k,\quad \forall\,k=1,2,\ldots\]
The pointwise convergence $\widehat{z}^{i_l}\to \widehat{y}$ 
implies the pointwise convergence $\z^{i_l}\to \z$, i.e.,
\begin{equation}
\label{pointwise1}
\lim_{l\to\infty}\,z^{i_l}_{\alpha\beta}\,\to\,z_{\alpha\beta}\qquad \forall\,\alpha\in\N^{n},\beta\in\N^p.
\end{equation}
Next, let $s\in\N$ be fixed.
From the pointwise convergence (\ref{pointwise1}) we deduce that
\[\lim_{l\to\infty}\M_{s}(z^{i_l})\,=\,\M_s(\z)\,\succeq0.\]
Similarly
\[\lim_{l\to\infty}\M_{s}(h_j\,\z^{i_l})\,=\,\M_{s}(h_j\,\z)\,
\succeq0,\quad j=1,\ldots,t.\]
As $s$ was arbitrary, we obtain
\begin{equation}
\label{put2}
\M_s(\y)\succeq0;\quad \M_s(h_j\,\z)\succeq0,\quad j=1,\ldots,t; \quad s=0,1,2,\ldots,
\end{equation}
which by Theorem \ref{thput}
implies that $\z$ is the sequence of moments of some finite measure $\mu^*$ with support contained in $\K$.
Moreover, the pointwise convergence (\ref{pointwise1}) also implies that
 \begin{equation}
 \label{last}
 \int_\Y\y^\beta\,d\varphi(\y)\,=\,\gamma_\beta\,=\,
 \lim_{l\to\infty}z^{i_l}_{0\beta}\,=\,z_{0\beta}\,=\,\int_\K \y^\beta\,d\mu^*,\quad\forall\beta\in\N^p.
 \end{equation}
 As measures on compacts sets are determinate, (\ref{last}) implies
 that the marginal of $\mu^*$ on $\R^p$ is the probability measure $\varphi$, and so $\mu^*$ is feasible for $\P$.
Finally, combining the pointwise convergence (\ref{pointwise1}) with (\ref{value}) yields
\[\rho\,\geq\,\lim_{l\to\infty}\,\rho_{i_l}\,= \,\lim_{i\to\infty}\,L_{\z^{i_l}}(f)\,=\,L_\z(f)\,=\,
\int_\K f\,d\mu^*,\]
which in turn yields that $\mu^*$ is an optimal solution of $\P$. And so $\rho_{i_l}\to\rho$ as $l\to\infty$.
As the sequence $(\rho_{i})$ is monotone this yields the desired result (a).

(b) Next, let $\alpha\in\N^n$ and $\beta\in\N^p$ be fixed, arbitrary.
From (\ref{pointwise1}), we have:
\[\lim_{l\to\infty}\,z^{i_l}_{\alpha\beta}\,=\,z_{\alpha\beta}\,=\,\int_\K\bx^\alpha\,\y^\beta\,d\mu^*,\]
and by Theorem \ref{th1}(c)
\[\lim_{l\to\infty}\,z^{i_l}_{\alpha\beta}\,=\,\int_\K\bx^\alpha\,\y^\beta\,d\mu^*\,=\,
\int_\Y\y^\beta\,g(\y)^\alpha\,d\varphi(\y),\]
and as the converging subsequence was arbitrary, the above convergence holds for the whole sequence
$(z^i_{\alpha\beta})$. $\qed$

\end{document}